# MAFIA: A THEORETICAL STUDY OF PLAYERS AND COALITIONS IN A PARTIAL INFORMATION ENVIRONMENT

By Mark Braverman,[1] Omid Etesami[2] and Elchanan Mossel[3]

*University of Toronto, University of California, Berkeley and University of California, Berkeley*

In this paper, we study a game called "Mafia," in which different players have different types of information, communication and functionality. The players communicate and function in a way that resembles some real-life situations. We consider two types of operations. First, there are operations that follow an open democratic discussion. Second, some subgroups of players who may have different interests make decisions based on their own group interest. A key ingredient here is that the identity of each subgroup is known only to the members of that group.

In this paper, we are interested in the best strategies for the different groups in such scenarios and in evaluating their relative power. The main focus of the paper is the question: How large and strong should a subgroup be in order to dominate the game?

The concrete model studied here is based on the popular game "Mafia." In this game, there are three groups of players: Mafia, detectives and ordinary citizens. Initially, each player is given only his/her own identity, except the mafia, who are given the identities of all mafia members. At each "open" round, a vote is made to determine which player to eliminate. Additionally, there are collective decisions made by the mafia where they decide to eliminate a citizen. Finally, each detective accumulates data on the mafia/citizen status of players. The citizens win if they eliminate all mafia members. Otherwise, the mafia wins.

We first find a randomized strategy that is optimal in the absence of detectives. This leads to a stochastic asymptotic analysis where it is

Received May 2007; revised June 2007.
[1]Supported in part by an NSERC CGS scholarship. Part of the work was done while interning with the Theory Group at Microsoft Research.
[2]Supported by BSF Grant 2004105.
[3]Supported by an Alfred Sloan fellowship in Mathematics, BSF Grant 2004105 and NSF Grants DMS-05-28488 and DMS-05-48249 (CAREER) and by DOD ONR: Combinatorial Statistics on Trees and Networks (N0014-07-1-05-06).

*AMS 2000 subject classifications.* 91A18, 91A28, 60J20.

*Key words and phrases.* Mafia, extensive game, optimal strategy, martingale.

shown that the two groups have comparable probabilities of winning exactly when the total population size is $R$ and the mafia size is of order $\sqrt{R}$.

We then show that even a single detective changes the qualitative behavior of the game dramatically. Here, the mafia and citizens have comparable winning probabilities only for a mafia size *linear* in $R$. Finally, we provide a summary of simulations complementing the theoretical results obtained in the paper.

## 1. Introduction.

*Motivation*: *quantitative analysis of extensive games.* In many real-life games, different players and coalitions have different information and actions available to them. Familiar examples include: workers, managers and stockholders of a company; students, teachers and management of a school; soldiers, citizens and insurgents in a war; or citizens, the mafia and the police in a certain city.

All of these games share similar features. For example, each player may belong to one or more coalitions. A worker/manager may also be a stockholder. A soldier/insurgent may be a citizen of the city where fighting takes place. And a citizen in a city may be a policeman or a mafia member.

Another common feature is that different groups make decisions in different ways and take actions of different types. In a company, a worker may influence the company's future by his/her direct actions at the company, while as a manager/shareholder, he/she may participate in various types of votes that will determine the future of the company.

The final common feature is the existence of different ways in which different players accumulate information. In particular, in all of the above examples, there is some public information that is available to all, as well as some other types of information that are available only to specific players/coalitions.

*The "Mafia" game.* A popular game exhibiting these different kinds of interactions is "Mafia," which takes place in an imaginary city. The mafia is trying to destroy this city. The mafia, citizens and detectives all have different information and available actions.

In this paper, we analyze this game with a particular focus on the relative power of the different groups. As we will see, the different players have immensely different powers: for a large population of size $R$, a mafia of order $\sqrt{R}$ already has a reasonable chance of winning—and a larger mafia will surely win.

Interestingly, as soon as there is one detective present, the game becomes fair only when the mafia consists of a linear fraction of the total population.

The fact that different kinds of information and actions yield a great variation in relative power is clearly true in many other settings. Our investigation here is an initial step toward understanding this phenomenon. In the concluding section, we discuss some more general insights resulting from our analysis.

We note, in particular, that previous research on partial information and group games is mostly concerned with general definitions and abstract results in the context of extensive games. This line of research has not resulted in much quantitative analysis (see, e.g., [3, 4]). Here, we focus on a very restricted setting, but obtain very precise results on the relative power of different groups.

In particular, the following phenomena should be valid in further generality (see also Section 7):

- In cases where there exists a distinguished group of size $M$ that has complete information and acts at all rounds playing against a group of players of size $R - M$ with no prior information that acts only at a fraction $\alpha$ of the rounds, it is expected that the two groups will have comparable winning probabilities if $M = R^\alpha$.
- As soon as the group of players with no prior information acquires information at a speed comparable with the speed at which actions are taken, for the two groups to have comparable winning probabilities, it is required that $M$ and $R$ are of the same order.

*Outline.* The model studied in this paper is defined in Section 2. The optimal strategies for the game without detectives are derived in Section 3. In Section 4, we use martingale arguments in order to show that in the game without detectives, the citizens and mafia have comparable winning probabilities when the mafia is of a size that is of order square root of the size of the total population. In Section 5, we analyze the game with detectives and show the dramatic effect of the information collected by the detectives. In Section 6, we provide more refined information on the game without detectives using simulation studies. Some general insights and future research directions are discussed in Section 7.

**2. Mafia game: definition.** We study the following model of the game "Mafia."

2.1. *The different groups.*

- There are $R$ players. Each player is a *resident*.
- $M$ of the $R$ players are *mafia members*. A non-mafia-member will be called a *citizen*.

- $D$ of the $R$ players are *detectives*. The sets of detectives and mafia members are disjoint. In particular, $M + D \leq R$. All of the detectives are citizens.

At the beginning of the game, the following information is given to each player.

- Each player is given his/her own identity, that is, he/she is told whether he/she is a mafia member or a citizen and whether he/she is a detective or not.
- Each mafia member is given the identity of all other mafia members.

No other information is given.

2.2. *The different rounds.* The game consists of iterations of the following three rounds until the game terminates.

**Residents:** Here, all of the residents pick one player to eliminate by a plurality vote. Thus, each resident is supposed to choose one person he/she wishes to eliminate. The person who receives the most votes is then eliminated and his/her identity is revealed. The vote takes place after a discussion between all residents. In cases of a tie, the identity of the person to be eliminated is chosen uniformly at random from among all players who received the maximal number of votes.

**Mafia:** In this round, the mafia choose a citizen to eliminate. This decision is made without any information leaking to the other players. The only information that becomes publicly available is the identity of the player eliminated and whether he/she was a detective or not.

**Detectives:** In this round, each detective queries the mafia/citizen status of one player. This status is then revealed only to the detective. No public information is revealed as a result of this exchange.

Note that after round $t$, there are $R_t = R - 2t$ residents. We denote by $M_t$ the number of mafia members after round $t$ and by $D_t$ the number of detectives after round $t$.

2.3. *Objectives*: *termination rules.* The game has two possible outcomes: either the mafia wins or the citizens win.

- The citizens win if all mafia members have been eliminated and there are still citizens alive.
- The mafia wins if all citizens have been eliminated when there are still mafia members alive.

Note that the objectives of the mafia and citizens are group objectives. Single players do not care if they are dead or alive, as long as their group wins.

2.4. *How do groups decide*? Note that in the Detectives round, each detective makes his/her own choice. We require each detective to return his/her choice in time polynomial in $R = R_0$ (the initial number of players). However, in the two other rounds, a group makes a decision. We proceed with the formal definition of group decisions.

For the Mafia round, this is easy. Since all mafia members have exactly the same objective and the same information, we may assume that the same rational choice is based on the information that was revealed up to that point. We further require that the mafia choose its victim in time polynomial in $R$.

In the Residents round, the situation is more involved. In particular, our analysis requires the following assumption.

ASSUMPTION 1. We assume that the citizens have a way of performing the following:

1. In the analysis of the game without detectives, we assume that all residents can send a message to all other residents simultaneously. If the game has a trustable moderator, this can be achieved by the moderator collecting messages from all residents and then displaying all messages (if a player does not send a message, the value of her message is 0). It can also be performed by means of a protocol where each player writes his message on a pad and then all pads are displayed simultaneously (again, the value of an undisplayed message is 0). Finally, this can also be implemented using commitment schemes under cryptographic assumptions [5].
2. In the analysis of the game with detectives, we need to assume that residents can vote anonymously. Given a trusted moderator, this can be achieved by a ballot run by the moderator. Otherwise, this can be performed using cryptographic voting schemes [1]
3. In the analysis of the game with detectives, we need to assume that residents can securely exchange messages (where it is only known how many messages were sent and received by each player). Given a trusted moderator, this can be achieved by letting him/her carry the messages. Otherwise, this can be achieved via standard public key techniques.

2.5. *Liveness*: *ensuring the game progresses.* We need to further specify the game's protocol to ensure the *liveness* property. This property is a common requirement in distributed protocols and software reliability and informally says that the game cannot "stall" (see, e.g., [6, 7]).

By the requirements above, each Detective and Mafia round lasts a polynomial number of steps. We model the "discussion" during the Residents round using communication rounds between the players. We assume that there is an order on the players and that in each communication round, each player has an opportunity to communicate openly to everyone and also to send private messages to other players. The length of all messages is polynomial in $R$.

Under Assumption 1.2 above, we allow the players to conduct an anonymous vote in one communication round. Note that the anonymous vote is not binding by itself and that, ultimately, the player to be eliminated at the end of the round is determined in an open plurality vote.

We further require each Residents round to last a polynomial number of steps. In other words, there is a $c$ such that for each Residents round, a vote is performed after $p(R) = O(R^c)$ communication rounds between the players. After $p(R)$ rounds, each player is *required* to vote openly in a predetermined order on the next player to be eliminated. We further bound the amount of computational steps each player can undertake between the communication rounds by $O(R^c)$.

The conditions above ensure that each Residents round takes at most $O(R^{2c})$ steps and that the entire game terminates in polynomial time. In practice, the protocols analyzed here satisfy the requirement above as they only require a constant number of communication rounds with a linear amount of computation (with an overhead added depending on the cryptographic protocols used). The protocols rely on analyzing and controlling the flow of information, rather than on complicated communication schemes.

**3. The game without detectives: optimal strategies.** In this section, we demonstrate that the game without detectives has a simple optimal strategy for both sides.

3.1. *Citizens' optimal strategy.* The citizens' strategy is designed in such away that if all citizens follow it then they have a high probability of winning the game. More specifically, this strategy guarantees that a random player will be eliminated as long as there is a majority of citizens. The strategy is defined as follows.

- On day $t$, each resident $1 \leq s \leq R_t$ picks a random integer $n(s)$ between 0 and $R_t - 1$. The residents announce their number simultaneously. Recall that the residents can announce their numbers simultaneously by Assumption 1.
- Let $n = 1 + (\sum n(s) \bmod R_t)$. All residents are supposed to vote to eliminate player number $n$.

3.2. *Properties of citizens' strategy.* Note that following this strategy will result in eliminating a random player as long as the citizens form a majority. This follows since the number $n$ has the uniform distribution as long as there is at least one citizen and since, when there is a majority of citizens, the vote will result in the elimination of player number $n$.

Note, furthermore, that the above protocol relies on the assumption that the players can announce their numbers at the same time before they may observe other players' announced numbers since, otherwise, the mafia may vote for a number such that the sum corresponds to a citizen.

CLAIM 1. *The strategy above is optimal for the citizens.*

This claim follows since, for every possible strategy for the citizens, all mafia members may follow this strategy, pretending to be citizens, until the mafia has achieved a majority. Note that by doing so, in each Residents round, a random resident will be eliminated, as in the strategy above. Once the mafia has achieved a majority, it will win, regardless of the citizens' strategy.

### 3.3. Mafia's optimal response.

CLAIM 2. *Any strategy for the mafia where all mafia vote to eliminate resident number $n$ that has been selected during the Residents round is an optimal response to the citizens' strategy.*

The second claim follows since as long as the citizens have the majority, the actions of the mafia are irrelevant. Moreover, once the mafia have a majority and they eliminate a citizen at each round, the mafia will win the game.

**4. The game without detectives: stochastic analysis.** Given the optimal strategies described above, the analysis of the game with detectives reduces to the analysis of the following stochastic process. Suppose that after round $t$, there are $R_t$ residents, of which $M_t$ are mafia members. Note that at the Residents round, a mafia member is eliminated with probability $M_t/R_t$ and that at the Mafia round, no mafia members are eliminated.

DEFINITION 1. Let $w(R, M)$ denote the probability that the mafia wins the game without detectives when initially there are $M$ mafia members among the $R$ residents and the citizens play according to their optimal strategy.

The following theorem roughly states that when there are no detectives, the mafia and citizens have comparable chances to win when the mafia size $M$ is of order $\sqrt{R}$. Moreover, if $M$ is a large multiple of $\sqrt{R}$, then the chance that the mafia wins is close to 1 and if it is a small multiple of $\sqrt{R}$, then the chance that the mafia wins is close to 0.

THEOREM 1. *There exists functions $p:(0,\infty) \to (0,1)$ and $q:(0,\infty) \to (0,1)$ such that if $0 < \eta < \infty$, the number of residents $R$ is sufficiently large and the mafia size satisfies $M \in [\eta\sqrt{R}, \eta\sqrt{R}+1]$, then*

$$p(\eta) \leq w(R,M) \leq q(\eta).$$

*Furthermore,*

$$\lim_{\eta \to \infty} p(\eta) = 1$$

*and*

$$\lim_{\eta \to 0} q(\eta) = 0.$$

We prove the two parts of Theorem 1 in Claim 3 and Claim 4 below.

CLAIM 3. *For every constant $\eta > 0$, there exists a constant $q(\eta) < 1$ such that for large enough $n$, when the Mafia has $M \leq \eta\sqrt{R}$ members, the Mafia will win with probability at most $q(\eta)$. Moreover, we have $\lim_{\eta \to 0} q(\eta) = 0$.*

PROOF. Let $R_t$ and $M_t$ denote the numbers of residents and mafia members, respectively, at the beginning of day $t$. The sequence $X_t = M_t(M_t - 1)/R_t$ is a martingale:

$$E[X_{t+1}|X_t, X_{t-1}, \ldots, X_0]$$
$$= \frac{M_t}{R_t} \cdot \frac{(M_t - 1)(M_t - 2)}{R_t - 2} + \frac{R_t - M_t}{R_t} \cdot \frac{M_t(M_t - 1)}{R_t - 2} = X_t.$$

We stop the martingale at the stopping time $T$ when either:

- there is at most one mafia member or at most one citizen, or
- the number of residents is less than or equal to $n+1$ for $n = n(\eta) = \lceil 4 + 8\eta^2 \rceil$.

By the martingale stopping time theorem (e.g., see [2]), we have

$$E[X_T] = E[X_0] = \frac{M(M-1)}{R} \leq \eta^2.$$

We now consider two cases.

- The case where $M_T < R_T/2$. Note that in this case, either $R_T = n, n+1$ or $M_T \leq 1$ and $R_T > n+1$. In the first case, we can bound the probability that the citizens win from below by the probability that at all Residents rounds, a mafia member is eliminated, which is bounded below by

$$\left(\frac{1}{n+1}\right)^{n/2+1}.$$

It is easy to see that the probability that the citizens win in the second case is at least $1/(n+1)$.

- The case where $M_T \geq R_T/2$. Note that if this is the case, then
$$X_T \geq (M_T - 1)/2 > n/4 - 1.$$

However,
$$P[X_T \geq n/4 - 1] \leq \frac{E[X_T]}{n/4 - 1} \leq \frac{\eta^2}{n/4 - 1}.$$

Thus, $P[X_T \geq n/4 - 1] \leq 1/2$ if $n \geq 4 + 8\eta^2$.

This proves that for every $\eta > 0$, the probability that the citizens win is at least
$$\frac{1}{2}\left(\frac{1}{n+1}\right)^{n/2+1},$$
proving that $q(\eta) < 1$ for all $\eta$.

Next, we want to show that if $\eta$ is sufficiently small, then $q(\eta)$ is close to 0. To achieve this, let $\epsilon > 0$. Now, repeat the argument above where:

- We first choose $n = n(\epsilon)$ large enough so that if $M_T \leq 1$ and $R_T \geq n$, then the mafia wins with probability at most $\epsilon/2$. This can be done since the probability that a single mafia member will win when there are $n$ residents is at most
$$g(n) = \left(1 - \frac{1}{n}\right)\left(1 - \frac{1}{n-2}\right) \leq \cdots \leq \left(\frac{n}{n+1}\frac{n-1}{n}\cdots\right)^{1/2} \leq \left(\frac{2}{n+1}\right)^{1/2}.$$

Thus, it suffices to take $n + 1 = 8/\epsilon^2$.

- We now repeat the argument above, considering the following two cases:
$$M_T \leq 1;$$
$$M_T \geq 2.$$

In the first case, the mafia will win with probability at most $\epsilon/2$. In the second case, if $M_T \geq 2$, then $X_T \geq 2/(n+1)$. On the other hand,
$$P[X_T \geq 2/(n+1)] \leq (n+1)\eta^2/2.$$

Thus, if $\eta \leq \min(2/(n+1), \epsilon/2)$, then we obtain that the second case occurs with probability at most $\epsilon/2$.

We have thus shown that if $\eta$ is sufficiently small, then the probability that the mafia wins is at most $\epsilon$, proving that $\lim_{\eta \to 0} q(\eta) = 0$. The quantitative estimates we obtain here show that it suffices to take $\eta \leq \epsilon^2/8$ in order to ensure that the mafia wins with probability at most $\epsilon$. □

CLAIM 4. *For every constant $\eta > 0$, there exists a constant $p(\eta) > 0$ such that for large enough $R$, when the Mafia has $M \geq \eta\sqrt{R}$ members, the Mafia wins with probability at least $p(\eta)$. Moreover, $\lim_{\eta \to \infty} p(\eta) = 1$.*

PROOF. Consider the sequence
$$Y_t = Y(R_t, M_t) := \frac{M_t^2(M_t-1)^2}{R_t^2 - R_t M_t + c M_t^2(M_t-1)^2},$$
where $c > 0$ is an appropriately chosen small constant. For example, one can take $c = 1/100$.

CLAIM 5. *There exists a $k > 0$ such that whenever $k \leq M_t < R_t/2$, then $Y_t$ is a submartingale, that is,*
$$E[Y_{t+1}|Y_t, Y_{t-1}, \ldots, Y_0, M_t \geq k] \geq Y_t.$$

PROOF. Given $M_t = M \geq 2$ and $R_t = R \geq 1$, it holds that
$$E[Y_{t+1}|M_t = M, R_t = R]$$
$$= E[Y(R_{t+1}, M_{t+1})|M_t = M, R_t = R]$$
$$= \frac{M}{R} \cdot Y(R-2, M-1) + \frac{R-M}{R} \cdot Y(R-2, M).$$

We claim that this quantity is greater than $Y(R, M)$. Denote the denominator of $Y(R, M)$ by $D(R, M) = R^2 - RM + c(M-1)^2 M^2$. Note that if $R \geq M > 0$, then $D(R, M)$ is positive. Using `Mathematica`, one obtains (with $c = 1/100$) that
$$E[Y_{t+1}|(R_t, M_t) = (R, M)] - Y(R, M)$$
is given by
$$\frac{P(R, M)/100}{R \cdot D(R, M) \cdot D(R-2, M) \cdot D(R-2, M-1)},$$
where
$$\begin{aligned}P(R, M) =\ & 1600 R^2 M - 1600 R^3 M + 400 R^4 M - 2416 RM^2 - 384 R^2 M^2 \\ & + 2400 R^3 M^2 - 800 R^4 M^2 + 16 M^3 + 4456 RM^3 - 4076 R^2 M^3 \\ & + 100 R^3 M^3 + \mathbf{400 R^4 M^3} - 72 M^4 - 1448 RM^4 + 2844 R^2 M^4 \\ & - 1000 R^3 M^4 + 122 M^5 - 847 RM^5 + 64 R^2 M^5 + \mathbf{100 R^3 M^5} \\ & - 88 M^6 + 308 RM^6 - 36 R^2 M^6 + 12 M^7 - 66 RM^7 - \mathbf{12 R^2 M^7} \\ & + 16 M^8 + 12 RM^8 - 6 M^9 + \mathbf{RM^9}.\end{aligned}$$

We need to show that for sufficiently large $k$, $P(R, M)$ is always positive. For sufficiently large $M$ and $R$, the highlighted terms of $P(R, M)$ dominate its behavior in the following sense. For each highlighted term, all of the monomials preceding it have both their $R$ and $M$ degree bounded by that of the highlighted term and at least one of them is strictly smaller. In particular, for $R \geq M \geq k$, for sufficiently large $k$, it holds that

$$P(R,M) > 300R^4M^3 + 99R^3M^5 - 13R^2M^7 + \tfrac{9}{10}RM^9$$
$$> 81R^3M^5 - 16R^2M^7 + \tfrac{64}{81}RM^9 = RM^5 \cdot (9R - \tfrac{8}{9}M^2)^2 \geq 0. \quad \square$$

We now return to the proof of Claim 4. We stop $Y$ at the first time $T$ where either:

- at least half of the remaining residents are Mafia members, or
- $M_t \leq k$.

We have

$$E[Y_T] \geq E[Y_0] = \frac{M^2(M-1)^2}{R^2 - RM + cM^2(M-1)^2}$$
$$\geq \frac{1}{R^2/(M^2(M-1)^2) + c}$$
$$\geq \frac{1}{2/\eta^2 + c} = \frac{1}{c} - \frac{2}{c(2+c\eta^2)} := h(\eta).$$

Observe that we always have $0 \leq Y_T < 1/c$. Moreover, if $M_T \leq k$, then

$$Y_T < \frac{k^4}{R_T(R_T - M_T)} < \frac{k^4}{R_T}.$$

Letting

$$A_1 = \left[M_T = k \text{ and } R_T > \frac{2k^4}{h(\eta)}\right],$$

we have

$$h(\eta) \leq E[Y_T] < P[A_1] \cdot \frac{k^4}{R_T} + (1 - p[A_1]) \cdot \frac{1}{c}$$
$$< \frac{h(\eta)}{2} + (1 - p[A_1]) \cdot \frac{1}{c}.$$

Hence,

$$1 - p[A_1] \geq \frac{c \cdot h(\eta)}{2}.$$

Thus, with probability at least $\frac{ch(\eta)}{2}$, either:

- at least half of the remaining residents are Mafia members, in which case the Mafia wins; or

- $M_T = k$ and $R_T \leq \frac{2k^2}{h(\eta)}$. The quantities here depend on $\eta$, but not on $R_0$, and the Mafia has a positive probability $s(\eta) > 0$ of winning with these initial conditions.

Finally, we obtain

$$p(\eta) \geq (1 - P[A_1]) \cdot s(\eta) \geq \frac{c \cdot h(\eta) s(\eta)}{2} > 0.$$

In order to conclude the proof, we would like to show now that $p(\eta) \to 1$ as $\eta \to \infty$. Let $A_2$ be the event that $M_T = k$. Then, on the complement of $A_2$, the Mafia wins. Therefore, it suffices to show that $P[A_2] \to 0$ as $\eta \to \infty$.

If $M_T = k$, then $R_T \geq 2M_T = 2k$ and therefore

$$Y_T < \frac{k^2(k-1)^2}{k^2 + ck^2(k-1)^2} = \frac{1}{c} - d(k)$$

for some positive $d(k)$. We now conclude that

$$h(\eta) \leq E[Y_T] < P[A_2]\left(\frac{1}{c} - d(k)\right) + (1 - P[A_2])\frac{1}{c}.$$

Since $h(\eta) \to \frac{1}{c}$ as $\eta \to \infty$, it follows that $P[A_2] \to 0$, and therefore $p(\eta) \geq 1 - P[A_2] \to 1$. □

## 5. The game with detectives.

5.1. *Results.* In this section, we investigate the power of the detectives. We show that even a single detective suffices to change the qualitative behavior of the game. More formally we prove the following.

THEOREM 2.

- Consider the game with one detective and mafia of size $M = \eta R < R/49$. Then, for $R$ sufficiently large, the probability that the mafia wins, denoted $w(R, M, 1)$, satisfies $p(\eta, 1) \leq w(R, M, 1) \leq q(\eta, 1)$, where $0 < p(\eta, 1) < q(\eta, 1) < 1$ for all $\eta < 1/49$ and $q(\eta, 1) \to 0$ as $\eta \to 0$.
- Let $d \geq 1$ and consider the game with $d$ detectives and mafia of size $M = \eta R$, where $\eta < 1/2$. Then, for $R$ sufficiently large, the probability that the mafia wins, denoted $w(R, M, d)$, satisfies $w(R, M, d) \leq q(\eta, d)$, where, for each $\eta < 1/2$, it holds that $\lim_{d \to \infty} q(\eta, d) = 0$.

The theorem implies that even a single detective dramatically changes the citizens' team power: while in the game with no detectives, a mafia of size $R^{1/2+\epsilon}$ will surely win, as soon as there is one detective, the mafia will lose unless it is of size $\Omega(R)$.

The rest of the section is dedicated to proving Theorem 2. In Section 5.2, we find a strategy for citizens that shows the existence of $q(\eta, 1) < 1$, such that $q(\eta, 1) \to 0$ as $\eta \to 0$. This will be proven in Claim 9. In Section 5.4, we find a strategy for the mafia that shows the existence of $p(\eta, 1) > 0$ for all $\eta < 1/2$. In fact, Claim 10 shows that such a strategy exists for any number of detectives. Finally, in Section 5.5, Claim 11, we prove the second part of the theorem.

5.2. *The citizens' strategy.* The key to the citizens' strategy is using the information gathered by the detective in an optimal way. Somewhat surprisingly, it turns out that the crucial information collected by the detective is not the identity of mafia members, but the identity of citizens. Note that the "natural" life expectancy of the detective is about $R/4$, in which time he/she will have a chance to collect information about roughly half of the Mafia members (the life expectancy is, in fact, smaller since citizens die faster than mafia). Even assuming that this half of the Mafia is eliminated, the remaining citizens will have to deal with the second half after the detective is gone. From previous sections, we know that this is impossible if the mafia size is $R^{1/2+\epsilon}$ for $\epsilon > 0$.

Instead, what is crucial is to use the information collected by the detective in order to notify citizens about the identity of other citizens. In this way, the citizens can collaborate after the detective has been eliminated in order to control the Residents rounds in a way similar to the way the mafia controls the Mafia rounds.

The citizens' strategy is divided into two stages. The first stage is before the detective dies and the second is after. Note that the death of the detective confirms his/her identity. The citizens strategy is as follows.

- **Stage 1: the detective is still alive.** This lasts for $\sqrt{\eta}R$ rounds, or until the detective is eliminated.
  - **The detective** collects information about people at random.
  - **The other citizens**, during the day phase, vote at random to eliminate a person, as in the case with no detectives.
- **Stage 2: the transition.** If the detective does not survive to this stage, the citizens forfeit. In other words, if the detective dies before round $\sqrt{\eta}R$, the citizens give up. Otherwise, the detective compiles an ordered list of people $V$ (vigilantes) that he/she knows are citizens. He/she then encrypts and sends the list to each member of $V$. At this stage, we are supposed to have $|V| > |M|$.

  The detective then asks everyone to eliminate him/her (during the day phase). Once the detective is eliminated (and thus the members of $V$ learn that the messages they have received are genuine), the third stage of the game begins.

　　In the case where multiple people claim to be detectives, they will all be eliminated according to the order in which they made their declarations (this guarantees that no mafia members will want to declare that they are detectives).
- **Stage 3: the detective is dead.** This lasts until the Mafia is eliminated, or until $|V| \leq |M|$, in which case the citizens forfeit. During every day round, the next person $p$ to be eliminated is selected using a *secure anonymous vote* and then everyone (at least the citizens) vote to eliminate $p$.
  – **Members of $V$:** The highest ranking surviving member of $V$ randomly selects a person from outside $V$ and sends his number $k$ to all other members of $V$. The members of $V$ then all vote for $k$ in the secret vote.
  – **Other citizens:** All other citizens abstain in the secret vote.

5.3. *Stochastic analysis of the citizens' strategy.*

CLAIM 6.　*The probability that the detective survives until the second stage is at least $1 - p_1(\eta)$, where*

$$p_1(\eta) = \frac{2\sqrt{\eta}}{1-\eta} < \frac{1}{3}$$

*and $p_1(\eta) \to 0$ as $\eta \to 0$.*

PROOF.　During the first stage of the strategy, the detective is indistinguishable from the rest of the citizens and his/her chance of being eliminated before making $\sqrt{\eta}R$ queries is bounded by $p_1(\eta) = \frac{2\sqrt{\eta}}{1-\eta} < \frac{1}{3}$. It is obvious that $p_1(\eta) \to 0$ as $\eta \to 0$. □

CLAIM 7.　*With probability at least $1 - p_1(\eta) - p_2(\eta)$, the detective survives until stage 2 and queries at least $v(\eta)$ of the citizens alive up to stage 2, where:*

- $p_2(\eta) \leq 2/3$ *for all $\eta < 1/49$ and $p_2(\eta) \to 0$ as $\eta \to 0$;*
- *the function $v(\eta)$ satisfies*

(1) $$v(\eta) \geq \tfrac{5}{2}\eta R = \tfrac{5}{2}M$$

*for all $\eta < 1/49$ and*

(2) $$v(\eta) \geq \sqrt{\eta}R/2 = \frac{M}{2\sqrt{\eta}}$$

*for small values of $\eta$.*

PROOF. In order to prove the claim, we will assume that the detective writes down the indices of the people he/she is going to query ahead of time. During rounds when he/she is supposed to query the identity of a dead resident, he/she will not query at all. Let $v(\eta)$ denote the number of citizens queried from the list that are alive by stage 2. Then, in order to prove the claim, it suffices to show that, with probability $1 - p_2(\eta)$, at least $v(\eta)$ of the citizens alive are stage 2 are on the list of queried residents.

Note that there are at least $(\sqrt{\eta} - \eta)R$ and at most $\sqrt{\eta}R$ citizens on the detective's list. At most $2\sqrt{\eta}R$ out of $(1-\eta)R$ citizens are eliminated. The set of those eliminated is chosen independently of the ones to be queried. Hence, the expected number of citizens that have been eliminated and also on the querying list is bounded by

$$\text{(3)} \qquad \left(\frac{2\sqrt{\eta}}{1-\eta}\right)\sqrt{\eta}R = \frac{2\eta}{1-\eta}R \leq \frac{7}{3}\eta R.$$

The last inequality holds for $\eta \leq 1/7$. Hence, with probability $\geq 1/3$, at most $(7/2)\eta R$ citizens on the list are eliminated, which means that at least

$$V \geq (\sqrt{\eta} - \tfrac{9}{2}\eta)R \geq \tfrac{5}{2}\eta R$$

survive to be in $V$. The last inequality assumes that $\eta < 1/49$. Together, we obtain that, with probability at least $1/3 - p_1(\eta) > 0$, the game survives to the second stage and $v(\eta) \geq \tfrac{5}{2}M$.

For a small $\eta$, using (3), the probability that at least $(\sqrt{\eta}/2 - \eta)R$ of the citizens on the list are eliminated is bounded by

$$p_2(\eta) = \frac{7/3\eta}{\sqrt{\eta}/2 - \eta} \to 0$$

as $\eta \to 0$. Hence, with probability at least $1 - p_2(\eta)$, we have

$$v(\eta) \geq (\sqrt{\eta} - \eta)R - (\sqrt{\eta}/2 - \eta)R = \sqrt{\eta}/2.$$

Thus, with probability at least $1 - p_1(\eta) - p_2(\eta) \to 1$ as $\eta \to 0$, the game survives to the second stage and $v(\eta) \geq M/(2\sqrt{\eta})$. □

CLAIM 8. *Consider stage 3 of the game with $|V| \geq v(\eta)$. Then, for all $\eta < 1/49$, the probability that the citizens lose is at most $p_3(\eta)$, where $p_3(\eta) \leq 4/5$ and $p_3(\eta) \to 0$ as $\eta \to 0$.*

PROOF. We define the time when the third stage begins as $t = 0$. Consider the quantity

$$Z_t = Z(V_t, M_t) = \frac{M_t}{V_t + 1}.$$

Define a stopping time $T$ to be the first time when either:

- $V_T \leq M_T$, in which case the citizens lose;
- there are no citizens outside $V$ remaining alive, in which case the citizens win if $V_T > M_T$;
- the mafia is eliminated, in which case the citizens win.

Denote by $U_t$ the citizens who are not members of $V$. Thus, during the course of the second stage, we have $U_t > 0$, $M_t > 0$ and $V_t > M_t$. We verify that under these conditions, $Z_t$ is a supermartingale. For any $V_t, M_t, U_t$ satisfying the conditions, we have

$$E[Z_{t+1}|V_t, M_t, U_t] - Z(V_t, M_t)$$
$$= \frac{M_t}{U_t + M_t} \cdot \frac{V_t}{U_t + V_t} Z(V_t - 1, M_t - 1)$$
$$+ \frac{M_t}{U_t + M_t} \cdot \frac{U_t}{U_t + V_t} Z(V_t, M_t - 1)$$
$$+ \frac{U_t}{U_t + M_t} \cdot \frac{V_t}{U_t + V_t - 1} Z(V_t - 1, M_t)$$
$$+ \frac{U_t}{U_t + M_t} \cdot \frac{U_t - 1}{U_t + V_t - 1} Z(V_t, M_t) - Z(V_t, M_t)$$
$$= \frac{M_t((M_t - V_t)(U_t + V_t) + 1 - M_t)}{(U_t + M_t)(V_t + 1)(U_t + V_t)(U_t + V_t - 1)} \leq 0.$$

The last inequality holds by our assumptions in the definition of the stopping time. Thus,

$$E[Z_T] \leq Z_0.$$

Observe that if the citizens lose, then $Z_T \geq 1/2$. In either case, $Z_T \geq 0$. Hence,

$$P[\text{citizens lose}] \leq \frac{E[Z_T]}{1/2} \leq 2Z_0.$$

In particular, we have the following.

- By (1), we have $V_0 \geq \frac{5}{2} M_0$ for all $\eta < 1/49$ and therefore

$$P[\text{citizens lose}] \leq 2Z_0 < 2 \cdot \tfrac{2}{5} = \tfrac{4}{5}.$$

  Hence, $p_3(\eta) \leq \tfrac{4}{5}$ for all $\eta < 1/49$.
- By (2), for small $\eta$, we have $V_0 \geq M/(2\sqrt{\eta})$ and therefore

$$P[\text{citizens lose}] \leq 2Z_0 < 2 \cdot 2\sqrt{\eta} = 4\sqrt{\eta} \to 0 \qquad \text{as } \eta \to 0.$$

  Thus, $p_3(\eta, 1) \to 0$ as $\eta \to 0$.

□

Since $q(\eta,1) \geq (1-p_1(\eta)-p_2(\eta))(1-p_3(\eta))$, Claim 6, Claim 7 and Claim 8 imply the following.

CLAIM 9. *The strategy defined in Section 5.2 satisfies $q(\eta,1) < 1$ for all $\eta < 1/49$ and $q(\eta,1) \to 0$ as $\eta \to 0$.*

5.4. *The Mafia's strategy.*

CLAIM 10. *In the game with $d$ detectives and mafia of size $\eta R$ for large $R$, the probability that the mafia wins is at least $\frac{\eta^2}{72}(\frac{\eta}{8d})^d$.*

PROOF. The mafia's strategy will be to eliminate random citizens. With probability at least $(\eta/8d)^d$, all of the detectives are dead by time $t_0 = \frac{\eta}{4d}R$. Moreover, by time $t_0$, there are at least $3\eta R/4$ mafia members that are alive and whose identity was not queried by any of the detectives. Finally, by time $t_0$, the detectives have queried the identity of at most $\eta R/4$ of the citizens.

The proof would follow if we could show that, given the scenario with a mafia of size at least $3\eta R/4$, the number of citizens whose identity was queried at most $\eta R/4$ and no detectives, the probability that the mafia wins is at least $\eta^2/72$. Let $V$ denote the set of citizens whose identities were queried and who are alive at time $t_0$. Let $S$ denote the mafia members alive at time $t_0$ and let $W$ denote the remaining citizens.

We first note that during the Residents rounds, the probability that a mafia member is eliminated is the same as the probability that a citizen whose identity has not been queried is eliminated. In other words, no matter what strategy the citizens choose, they cannot do better than eliminating at random one of the residents alive in the set $S \cup W$.

Let $S = S_1 \cup S_2$, where the sets $S_1$ and $S_2$ are disjoint and $|S_1| = |V|$. Note that $|S_2| \geq \eta R/2$. Consider the following suboptimal strategy of the mafia where, during the night, they eliminate uniformly a member of $V \cup W \cup S_2$.

Let $T$ be the first round where all member of $S_1$ or all members of $V$ have been eliminated. Then, by symmetry, it follows that, with probability at least $1/2$, all of the citizens in $V$ have been eliminated at time $T$.

Let $W(T)$ denote the number of citizens alive at time $T$. Let $S_2(T)$ denote the number of $S_2$ alive at time $T$. If $W(T) = 0$, then the mafia clearly wins. Otherwise, by conditioning on the value of $x = W(T) + S_2(T)$, we obtain

$$E[S_2(T)|W(T)+S_2(T) = x] \geq \frac{\eta/2}{1+\eta/2}x = \frac{\eta}{2+\eta}x > \frac{\eta}{3}x,$$

hence

$$E\left[\frac{S_2(T)}{W(T)+S_2(T)}\Big|W(T)+S_2(T) > 0\right] > \frac{\eta}{3}.$$

Therefore,

$$P\left[S_2(T) \geq \max\left(1, \frac{\eta}{6}W(T)\right)\right] \geq \frac{\eta}{6}.$$

Since the probability that a mafia of size $\max(1, \eta w/6)$, wins against $w$ citizens is at least $\eta/6$ the proof follows. □

5.5. *Many detectives.* We now prove the second part of Theorem 2.

CLAIM 11. *Consider the game with a mafia of size $\eta = (1/2 - \delta)R$. Let $d$ be an integer greater than $4/\delta + 1$. Then, given that there are at least $d^2$ detectives, the probability that the citizens win is at least $1 - de^{-d}$.*

PROOF. The strategy of the citizens and detectives is defined as follows. Before the game, the residents are partitioned into $d$ sets of size at most $\lfloor \frac{\delta}{4}R \rfloor$.

Each detective queries the identity of all players in a randomly chosen set. If the detective succeeds to query all identities, then he/she reveals all mafia members in that set. The detective is then eliminated in order to verify his/her claims.

If the identity of all mafia members has been revealed (by round $\frac{\delta}{4}R$), the citizens will win by eliminating one mafia member at each round since they are still a majority at round $\frac{\delta}{4}R$.

Note that a specific detective will not query the identity of one of the $d$ groups if either the detective is eliminated by the mafia or the detective picked a different group. Since the two events are independent, the probability that a specific group is queried by a specific detective is at least

$$\frac{1}{d}P[\text{the detective survives}] \geq \frac{1}{d}\left(1 - \frac{(\delta/2)R}{R - \eta R}\right)$$
$$> \frac{1}{d}\left(1 - \frac{(\delta/2)R}{R/2}\right) > \frac{1}{2d} \geq \delta.$$

This is also true conditioning on the status of all other detectives. Therefore, the probability that the status of the mafia members in the set is not queried is at most

$$(1 - \delta)^{d^2} \leq (1 - 2/d)^{d^2} \leq \exp(-d).$$

Therefore, the probability that all mafia members are queried is at least $1 - de^{-d}$. The proof follows. □

REMARK 1. It is an interesting problem to study the optimal querying procedure as a function of the number of detectives $d$ and the mafia size $\eta R$. In particular, some alternatives to the strategy suggested here include:

- the detectives will query identities at random;
- the detectives will query according to some combinatorial design.

We believe that for a small mafia size $\eta R$ and a small number of detectives, querying at random results in high winning probabilities. For high values of $\eta$ and $d$, it seems like that combinatorial designs should work better.

5.6. *A strategy for citizens with no cryptographic assumptions.* In this section, we briefly outline a strategy for citizens with at least one detective without making any cryptographic assumptions. The strategy gives the citizens a positive probability $p(\eta, d) > 0$ of winning against a mafia of size $\eta R$. Unlike the previous strategies, the present strategy makes no assumptions concerning private communication or other cryptographic protocols.

In this strategy, the detective collects information until time $T$, when he/she knows the identities of more than half of the residents alive to be good citizens. After that, the detective publishes a list $V_T$ of the good citizens he/she knows and is then eliminated to verify the claim. Under our assumption, $|V_T| > R_T/2$.

The citizens in $V_T$ then attempt during the day rounds to eliminate everyone not in $V_T$. Since $V_T$ are a majority, they will succeed and after $R_T - |V_T|$ rounds, only members of $V_T$ will remain alive and the citizens will win.

It remains to bound from below the probability that the detective will succeed without getting eliminated. The detective will make the queries independently at random.

CLAIM 12. *Assuming that $\eta < 1/72$, a detective has a probability of at least $p(\eta, 1) > 1/108$ to identify a set $V_T$ as above without being eliminated.*

PROOF. With probability at least $1/12$, the detective survives until the round when there are less than $L = R/9$ residents remaining. For a given surviving resident, his chance of not being queried in any round is at most

$$\left(1 - \frac{1}{R}\right)\left(1 - \frac{1}{R-2}\right) \cdots \left(1 - \frac{1}{L}\right)$$

$$< e^{-1/R} e^{-1/(R-2)} \cdots e^{-1/L} < e^{-(\ln R - \ln L)/2} = \sqrt{\frac{L}{R}} < \sqrt{\frac{1}{9}} = \frac{1}{3}.$$

Hence, we expect at least $\frac{2}{3}L$ of the residents to have been queried. This means that, with probability at least $\frac{1}{9}$, at least $\frac{5}{8}L$ residents have been queried. Since at most $R/72 = L/8$ of the residents can be mafia, we conclude that in this case, at least half of all residents have been identified as citizens. □

We conclude with the following observations.

- As the number of detectives grows, so do the chances of the above strategy to succeed. In particular, for any $\varepsilon > 0$, there exists a $d$ such that $d$ detectives have a winning probability of at least $1 - \varepsilon$ against a mafia of size $(1/2 - \varepsilon)R$.
- Unlike the strategy in Section 5.2, there is no guarantee that $p(\eta, d) \to 1$ as $\eta \to 0$. In fact, a mafia of size $\Omega(\sqrt{n})$ has a positive probability of winning—by eliminating the detectives before they have a chance to reveal information.

**6. Simulation studies.** In this section, we briefly discuss some simulations complementing the theoretical picture we have derived so far. All of the experiments deal with the case of no detectives, for which we know the optimal strategies of both the mafia and the citizens.

In Figure 1(a), we calculate the winning probability of a mafia of size $M = \eta\sqrt{R}$ as a function of $\eta$. The figure was derived by repeating the game 10,000 times with $R = 10,000$. In Figure 1(b), we zoom in to Figure 1(a) for $\eta \leq 0.4$ and simulating the game 20,000 times. Note that for such values, the function is almost linear.

In Figure 2, we estimate the size of the $M(R)$ such that the probability that the mafia wins is exactly $1/2$. This is done by running the game 2,000 times with different sizes of $M$ and $R$. Note the excellent fit of this function with a function of the type $M = c\sqrt{R}$.

In Figure 3, we run three simulations of the game with $R = 10^6$ and $M = 10^3$ (so $\eta = 1$). Each row in the figure corresponds to one run. In each of the three drawings, we plot the value of the martingale $X(t) = M_t \cdot (M_t - 1)/R_t$ as a function of the round $t$. In the first column we draw this function for all 500,000 steps, in the second for the last 10,000 and in the last for the last 100.

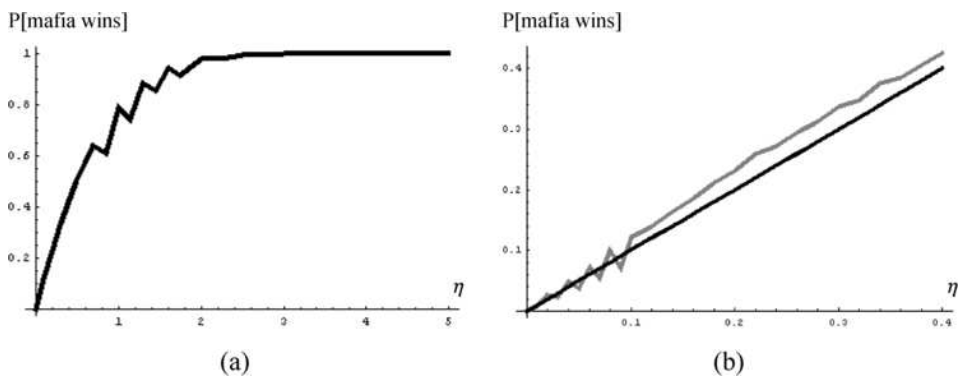

FIG. 1. (a) *The winning probability for mafia of size* $M = \eta\sqrt{R}$ *as a function of* $\eta$; (b) *The same winning probabilities for* $\eta \leq 0.4$—*the linear approximation of the probability is shown.*

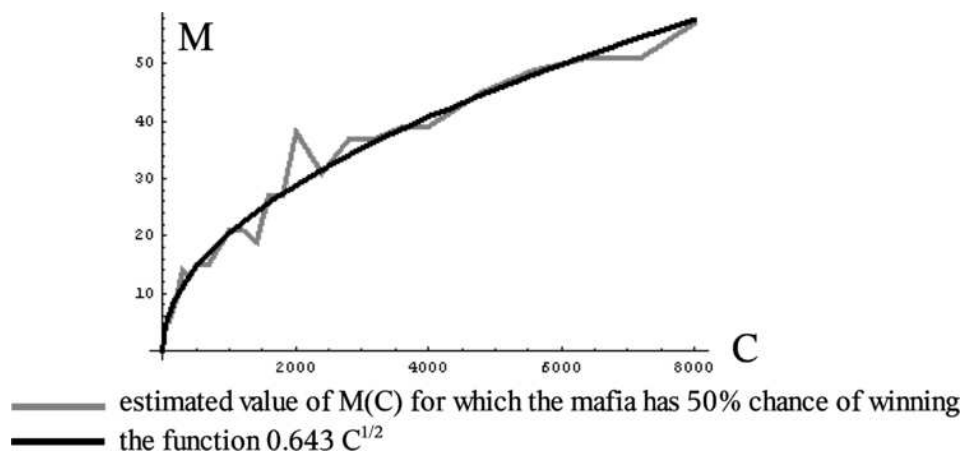

Fig. 2. *The size $M(R)$ for which the probability of the mafia to win is about $1/2$.*

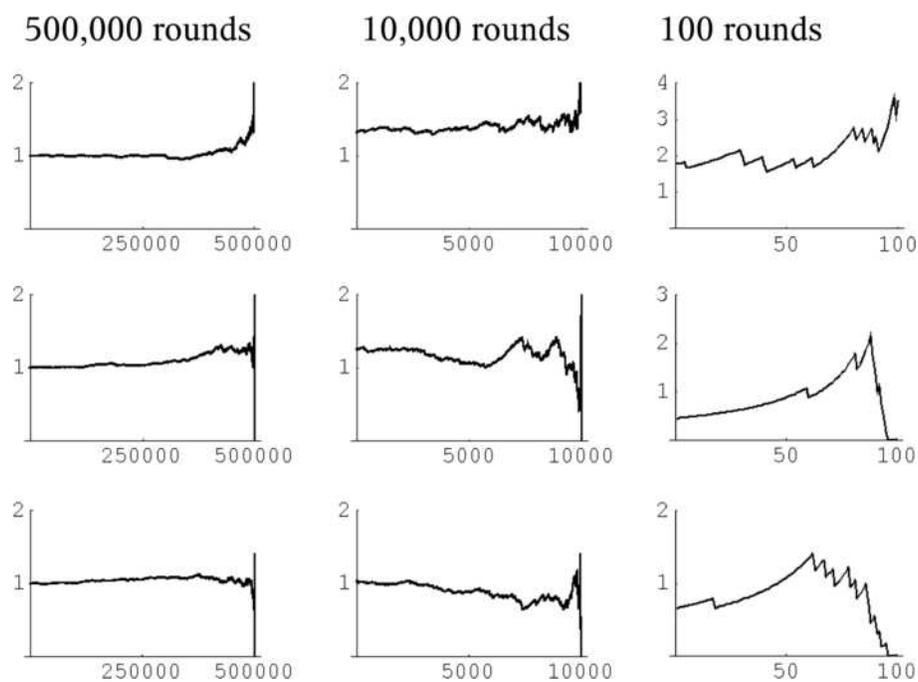

Fig. 3. *Three runs of the game where the function $X(t)$ of round $t$ is $X(t) = M_t \cdot (M_t - 1)/R_t$. The value of $X(t)$ is drawn for $t \leq 500{,}000$. The second and third columns zoom in on the last $10{,}000$ and $100$ rounds, respectively.*

From the figure, one sees that the function $\eta(t)$ is mostly "deterministic" at the beginning of the game and becomes "random" only at the game's end. While we have not stated such results formally, this can be shown using standard concentration results.

**7. Discussion and open problems.** Our results exhibit interesting trade-offs in a group game between information, actions and sizes of groups. They also raise questions about the significance of secure protocols in such games. We believe that these problems should be further investigated. In particular:

- In the mafia game without detectives, the mafia acted at all rounds, while the citizens acted at only half of the rounds.

  CONJECTURE 1. *Consider a variant of the mafia game without detectives, where each $r$ rounds are partitioned into $d$ day rounds and $n = r - d$ night rounds. The two groups then have comparable winning probabilities if $M = R^{d/r}$.*

  Note that our results correspond to the special case where $r = 2, d = 1$. In fact, one would expect such a phenomenon to be more general. It should hold when a group with complete information plays against a group where each individual has very little information and where the partial information group takes action at $d/r$ of the rounds.
- Our results show that once the partial information group can collect information at a linear speed, it stands a chance against even a complete information group of comparable size. It is interesting to study how general this phenomenon is.
- It would be interesting to see whether the strategy from Section 5.6 can be improved to guarantee success against a sublinear mafia with probability tending to 1 (as is the case *with* the cryptographic assumptions). If there is no such strategy, proving that this is the case appears to be very hard without putting strict restrictions on the type of messages allowed to be passed.

**Acknowledgments.** E. M. would like to thank Mohammad Mahdian, Vahab S. Mirrokni and Dan Romik for interesting discussions. This paper was conceived at the Miller Institute Annual Symposium, 2004, when E. M. was moderating a large Mafia game. E. M. is grateful to the Miller Institute for the opportunity to reflect on the mathematics of Mafia and other games.

M. Braverman
Department of Computer Science
University of Toronto
Toronto, Ontario
Canada M59 3G4
E-mail: mbraverm@cs.toronto.edu

O. Etesami
Department of Computer Science
University of California, Berkeley
Berkeley, California 94720-1776
USA
E-mail: etesami@cs.berkeley.edu

E. Mossel
Department of Statistics
University of California, Berkeley
Berkeley, California 94720-3860
USA
E-mail: mossel@stat.berkeley.edu